\documentclass[11pt]{article}
\usepackage{amscd,amsmath, amssymb, fancyhdr, color}

\newcommand{\version}{version 2.0,\ \  November 14, 2011}

\newcommand{\al}{\alpha}
\newcommand{\be}{\beta}

\newcommand{\RR}{\mathbb{R}}

\numberwithin{equation}{section}

\def\eqref#1{(\ref{#1})}

\newcommand{\ra}{{\:\longrightarrow\:}}

\newcommand{\C}{{\mathbb C}}
\newcommand{\R}{{\mathbb R}}

\def\1{\sqrt{-1}\:}
\newcommand{\restrict}[1]{{\left|_{{\phantom{|}\!\!}_{#1}}\right.}}
\newcommand{\cntrct}                
{\hspace{2pt}\raisebox{1pt}{\text{$\lrcorner$}}\hspace{2pt}}
\newcommand{\arrow}{{\:\longrightarrow\:}}

\renewcommand{\bar}{\overline}
\renewcommand{\phi}{\varphi}
\renewcommand{\epsilon}{\varepsilon}
\renewcommand{\geq}{\geqslant}

\newcommand{\im}{\operatorname{im}}

\newcommand{\const}{\operatorname{\text{\sf const}}}
\newcommand{\Vol}{\operatorname{Vol}}

\newcommand{\Map}{\operatorname{Map}}

\newcounter{Mycounter}[section]
\newcounter{lemma}[section]
\renewcommand{\thelemma}{{Lemma \thesection.\arabic{lemma}}}
\newcommand{\lemma}{%
     \setcounter{lemma}{\value{Mycounter}}
     \refstepcounter{lemma}
     \stepcounter{Mycounter}
     {\noindent \bf \thelemma:\ }}

\newcounter{claim}[section]
\renewcommand{\theclaim}{{Claim \thesection.\arabic{claim}}}
\newcommand{\claim}{%
     \setcounter{claim}{\value{Mycounter}}
     \refstepcounter{claim}
     \stepcounter{Mycounter}
     {\noindent \bf \theclaim:\ }}

\newcounter{sublemma}[section]
\newcounter{corollary}[section]
\renewcommand{\thecorollary}{{Corollary \thesection.\arabic{corollary}}}
\newcommand{\corollary}{%
     \setcounter{corollary}{\value{Mycounter}}
     \refstepcounter{corollary}
     \stepcounter{Mycounter}
     {\noindent \bf \thecorollary:\ }}

\newcounter{theorem}[section]
\renewcommand{\thetheorem}{{Theorem \thesection.\arabic{theorem}}}
\newcommand{\theorem}{%
     \setcounter{theorem}{\value{Mycounter}}
     \refstepcounter{theorem}
     \stepcounter{Mycounter}
     {\noindent \bf \thetheorem:\ }}

\newcounter{conjecture}[section]
\newcounter{proposition}[section]
\newcounter{definition}[section]
\renewcommand{\thedefinition}
       {{Definition~\thesection.\arabic{definition}}}
\newcommand{\definition}{%
     \setcounter{definition}{\value{Mycounter}}
     \refstepcounter{definition}
     \stepcounter{Mycounter}
     {\noindent \bf \thedefinition:\ }}

\newcounter{example}[section]
\renewcommand{\theexample}{{Example \thesection.\arabic{example}}}
\newcommand{\example}{%
     \setcounter{example}{\value{Mycounter}}
     \refstepcounter{example}
     \stepcounter{Mycounter}
     {\noindent \bf \theexample:\ }}

\newcounter{remark}[section]
\renewcommand{\theremark}{{Remark \thesection.\arabic{remark}}}
\newcommand{\remark}{%
     \setcounter{remark}{\value{Mycounter}}
     \refstepcounter{remark}
     \stepcounter{Mycounter}
     {\noindent \bf \theremark:\ }}

\newcounter{problem}[section]
\newcounter{question}[section]
\renewcommand{\thequestion}{{Question \thesection.\arabic{question}}}
\newcommand{\question}{%
     \setcounter{question}{\value{Mycounter}}
     \refstepcounter{question}
     \stepcounter{Mycounter}
     {\noindent \bf \thequestion:\ }}

\makeatletter

\@addtoreset{equation}{section}
\@addtoreset{footnote}{section}
\makeatother

\def\blacksquare{\hbox{\vrule width 5pt height 5pt depth 0pt}}
\def\endproof{\blacksquare}

\newcommand{\samethanks}[1][\value{footnote}]{\footnotemark[#1]}

\def\ddbar{\partial \overline{\partial}}

\begin{document}
\begin{center}
{\LARGE\bf
Blow-ups of locally conformally\\[3mm] K\"ahler manifolds}\\[3mm]
{\large 
Liviu Ornea\footnote{Partially supported by CNCS UEFISCDI, project
number PN-II-ID-PCE-2011-3-0118.}, 
Misha
Verbitsky\footnote{Partially supported by RFBR grant
10-01-93113-NCNIL-a, RFBR grant 09-01-00242-
a, AG Laboratory SU-HSE, RF government grant, ag. 11.G34.31.0023, and
Science Foundation
of the SU-HSE award No. 10-09-0015.}, and 
Victor Vuletescu\samethanks[1]}\\[3mm]

{\bf Keywords:} Locally
conformally K\"ahler manifold,
locally trivial bundle, blow-up.

{\bf 2000 Mathematics Subject
Classification:} { 53C55.}\\[4mm]

\end{center}

{\small
\hspace{0.15\linewidth}
\begin{minipage}[t]{0.7\linewidth}
{\bf Abstract} \\ A locally conformally K\"ahler (LCK) manifold is a manifold which is covered by a K\"ahler manifold, with the deck transform group acting by homotheties. We show that the blow-up of a compact LCK manifold along a complex submanifold admits an LCK structure if and only if this submanifold is globally conformally K\"ahler. We also prove that a twistor space (of a compact 4-manifold, a quaternion-K\"ahler manifold or a Riemannian manifold)  cannot admit an LCK metric, unless it is K\"ahler.
\end{minipage}
}
\tableofcontents
\section{Introduction}

\subsection{Bimeromorphic maps and locally conformally K\"ahler structures}

A {\bf locally conformally K\"ahler} (LCK) manifold is a 
complex manifold $M$, $\dim_\C M >1$, admitting
a K\"ahler covering $(\tilde M, \tilde \omega)$, with
the deck transform group acting on $(\tilde M, \tilde \omega)$
by holomorphic homotheties. Unless otherwise stated, we shall consider only compact LCK manifolds.

In the present paper we are interested in
the birational (or, more precisely, bimeromorphic)
geometry of LCK manifolds. 

An obvious question arises immediately.

\hfill

\question
Let $X\subset M$ be a complex subvariety of an LCK
manifold, and $M_1 \arrow M$ a blowup of $M$ in $X$.
Would $M_1$ also admit an LCK structure?

\hfill

When $X$ is a point, the question is answered in
affirmative by Tricerri \cite{_Tricerri_} and Vuletescu
\cite{vuli}. When $\dim X >0$, the answer is not
immediate. To state it properly, we recall the notion
of a {\bf weight bundle} of an LCK manifold. 
Let $(\tilde M, \tilde \omega)$ be the K\"ahler covering of an LCK
manifold $M$, and $\pi_1(M) \arrow \Map(\tilde M, \tilde M)$
the deck transform map. Since 
$\rho^*(\gamma) \tilde\omega = \const \cdot \tilde \omega$,
this  constant defines a character $\pi_1(M) \stackrel \chi \arrow \R^{>0}$,
with $\chi(\gamma) := \displaystyle\frac {\rho^*(\gamma)\tilde \omega}{\tilde \omega}$.

\hfill

\definition
Let $L$ be the 1-dimensional 
local system on $M$ with monodromy defined by the
character $\chi$. We think of $L$ as of a real bundle
with a flat connection. This bundle is called
{\bf the weight bundle} of $M$.

\hfill

One may think of the K\"ahler form 
$\tilde \omega$ as of an $L$-valued
differential form on $M$. This form
is closed, positive, and of type (1,1).
Therefore, for any smooth complex subvariety
$Z\subset M$ such that $L\restrict Z$ is 
a trivial local system, $Z$ is K\"ahler. 

\hfill

The following two theorems describe how
the LCK property behaves under blow-ups.

\hfill

\theorem \label{_blow_ups_intro_Theorem_}
Let $Z\subset M$ be a compact complex
submanifold of an LCK manifold, and 
$M_1$ the blow-up of $M$ with center in $Z$.
{
If the
restriction $L\restrict Z$ of the weight bundle
is trivial as a local system then $M_1$ admits an LCK metric.}

\hfill

{\bf Proof:} See \ref{_blow-up_LCK_Corollary_}. \endproof

\hfill

A similar question about blow-downs is also answered.

\hfill

\theorem \label{_blow_downs_intro_Theorem_}
Let $D\subset M_1$ be an exceptional divisor
on an LCK manifold, an $M$ the complex variety
obtained as a contraction of $D$. Then the
restriction $L\restrict D$ of the weight bundle
to $D$ is trivial.

\hfill

{\bf Proof:} \ref{_main_blow-down_Theorem_}. \endproof

\hfill

This result is quite unexpected, and leads to the following
theorem about a special class of LCK manifold called
{\em Vaisman manifolds} (Section \ref{_main_stateme_Section_}). 

\hfill

\claim 
Let $M$ be a Vaisman manifold.
Then any bimeromorphic contraction $M \arrow M'$
is trivial. Moreover, for any positive-dimensional
submanifold $Z\subset M$, its blow-up $M_1$
does not admit an LCK structure.

\hfill

{\bf Proof:} \ref{_Vaisman_blow_Corollary_} \endproof

\subsection{Positive currents on LCK manifolds}\label{cur}

The proofs of \ref{_blow_downs_intro_Theorem_} 
and \ref{_blow_ups_intro_Theorem_} are purely topological.
However, they were originally obtained using a less elementary 
argument involving positive currents. 

We state this argument here, omitting minor details
of the proof, because we think that this line of thought could be fruitful in other contexts too; 
for more information and missing details, 
the reader is referred to  \cite{_Demailly:ecole_}, 
\cite{_Demailly_Paun_} and \cite{_Demailly:Reg_}.

A {\bf current} is a form taking values in distributions.
The space of $(p,q)$-currents on $M$ is denoted by $D^{p,q}(M)$.
A {\bf strongly positive current}\footnote{In the present
paper, we shall often omit ``strongly'', because we are 
only interested in strong positivity.}
 is a linear combination
\[
\sum_I \alpha_I (z\wedge\bar z)_I
\]
where $\alpha_I$ are positive, measurable functions,
and the sum is taken over all multi-indices $I$.
An integration current of a
closed complex subvariety is a strongly
positive current. 

It is easy to define the de Rham differential on 
currents, and check that its cohomology coincide
with the de Rham cohomology of the manifold.

Currents are naturally dual to differential forms with compact
support. This allows one to define an integration 
(pushforward) map of currents, dual to the
pullback of differential forms. This map is
denoted by $\pi_*$, where $\pi:\; M \arrow N$
is a proper morphism of smooth manifolds.

Now, let $\pi:\; M \arrow N$ be a blow-up of a subvariety
 $Z\subset N$ of codimension $k$, and $\omega$ a K\"ahler
 form on $M$. Then $(\pi_*\omega)^k$ has
 a singular part which is proportional to
 the integration current of $Z$.
 
 This follows from the Siu's decomposition of
positive currents (\cite{_Demailly:ecole_}). 
Demailly's results on  intersection theory
of positive currents (\cite{_Demailly:Reg_})
are used to multiply the currents, and the rest follows
because the Lelong numbers of $\pi_*\omega$
 along $Z$ are non-zero.
 
Applying this argument to a birational
contraction $M \stackrel \phi \arrow M'$ of an LCK manifold $M$, 
and denoting by $\tilde M \stackrel {\tilde \phi} 
\arrow \tilde M'$ the corresponding map of coverings, 
we obtain a closed, positive current 
$\xi:= \tilde \phi_* \tilde \omega$ on  
$\tilde M'$, with the deck transform map $\rho$ acting on $\xi$ by homotheties.
Then $\rho$ would also act by homotheties on the
current $\xi^k$, $k =\dim Z$, where $Z$ is 
the exceptional set of $\tilde \phi$.

Applying the above result to decompose
$\xi^k$ onto its absolutely continuous
and singular part, we obtain that the
current of integration $[Z]$ of $Z$ is mapped
to $\const [Z]$ by the deck transform action.
Since the current 
of integration of $Z$ is mapped by the deck transform 
to the  current of integration of $\tilde \phi(Z)=Z$,
the constant $\const$ is trivial; this implies
that $\pi(Z)\subset M'$ is K\"ahler, with the
K\"ahler metric obtained in the usual way from
$\tilde \omega$.

\subsection{Fujiki class C and LCK geometry}

A compact complex variety $X$ is said to belong to
{\bf Fujiki class C} if $X$ is bimeromorphic to a 
K\"ahler manifold. The Fujiki class C manifolds are closed under
many natural operations, such as taking a subvariety, or
the moduli of subvarieties, and play important role
in K\"ahler geometry. 

This notion has a straightforward LCK analogue.

\hfill

\definition
Let $M$ be a compact complex variety. It is called
{ \bf a locally conformally class C variety} if it is 
bimeromorphic to an LCK manifold.

\hfill

The importance of the Fujiki class C notion
was emphasized by a more recent work of Demailly
and P\u aun \cite{_Demailly_Paun_}, 
who characterized class C manifolds in terms
of positive currents. Recall that a {\bf K\"ahler current}
is a positive, closed (1,1)-current $\phi$ on a complex
manifold $M$ which satisfies $\phi \geq \omega$ for some
Hermitian form $\omega$ on $M$. 

Demailly and P\u aun have proven that a compact complex
 manifold $M$ belongs to class C if and only if it admits
a positive K\"ahler current.

For an LCK manifold, an analogue of
a K\"ahler current is provided by the following notion (motivated by \ref{def2}).

\hfill

\definition
Let $M$ be a compact complex manifold,
$\theta$ a closed real 1-form on $M$, $\Xi$ a
positive, real (1,1)-current satisfying 
$d\Xi=\theta\wedge\Xi$ and $\Xi \geq \omega$
for some Hermitian form $\omega$ on $M$.
Then $\Xi$ is called {\bf an LCK current}.

\hfill

It would be interesting to know if an LCK-analogue
of the Demailly-P\u aun theorem is true.

\hfill

\question
Let $M$ be a complex compact manifold. Determine
whether the following conditions are equivalent.
\begin{description}
\item[(i)] $M$ belongs to 
locally conformally class C.
\item[(ii)] $M$ admits an LCK current.
\end{description}

\section{Blow-ups and blow-downs of LCK manifolds}
\label{_main_stateme_Section_}

We start by repeating (in a more technical fashion)
the definition of an LCK manifold given in the introduction.
Please see \cite{drag} for more details and several other 
versions of the same definition, all of them equivalent.

\hfill

\definition A {\bf locally conformally K\"ahler} (LCK) manifold is a
complex manifold
$X$ covered by a system of open subsets $U_\al$ endowed with {\em
local} K\"ahler metrics $g_\al$, conformal on overlaps $U_\al\cap U_\be$:
$g_\al=c_{\al\be}g_\be$.

\hfill

Note that, in complex dimension at least $2$, as
we always assume, $c_{\al\be}$ are positive constants. Moreover, they
obviously satisfy the cocycle condition. Interpreted in cohomology, the
cocycle
$\{c_{\al\be}\}$ determines a closed one-form $\theta$, called {\bf the Lee
form}. Hence, locally $\theta=df_\al$. It is easily seen that
$e^{-f_\al}g_\al=e^{-f_\be}g_\be$ on  $U_\al\cap U_\be$, and thus
determine a {\em global} metric $g$ which is conformal on each $U_\al$
with a K\"ahler metric.  One obtains the following equivalent:

\hfill

\definition\label{def2} A Hermitian manifold $M$ is LCK if its fundamental two-form
$\omega$ satisfies:
\begin{equation}
d\omega=\theta\wedge\omega, \quad \quad d\theta=0.
\end{equation}
for a {\em closed} one-form $\theta$.

\hfill

If $\theta$ is exact
then $M$ is called {\bf globally conformally
K\"ahler} (GCK).

As we work with compact manifolds and, in general, the topology of compact
K\"ahler manifolds is very different from the one of compact LCK
manifolds, we always assume $\theta\neq 0$ on $X$.

\hfill

Let $\Gamma \ra\tilde M\stackrel{\pi}\ra M$ be the
universal cover of $M$ with deck group $\Gamma$. As $\pi^*\theta$ is exact
on $\tilde M$,
$\pi^*\omega$ is globally conformal with a K\"ahler metric $\tilde
\omega$. Moreover, $\Gamma$ acts by
holomorphic homotheties with respect to $\tilde \omega$. This defines a
character
\begin{equation}\label{chi}
\chi:\Gamma\ra \R^{>0}, \quad \gamma^*\tilde\omega=\chi(\gamma)\tilde\omega.
\end{equation}
It can be shown that this property is indeed an equivalent definition of
LCK manifolds, see \cite{ov11}.

Clearly, a LCK manifold $M$ is globally conformally
K\"ahler if and only if $\Gamma$ acts trivially on
$\tilde \omega$ ({\em i.e.} $\im \chi=\{1\}$).

\hfill

A particular class of LCK manifolds are the {\bf Vaisman manifolds}. They
are LCK manifolds with the Lee form parallel with respect to the
Levi-Civita connection of the LCK metric. The compact ones are mapping
tori over the circle with Sasakian fibre, see \cite{ov03}. The typical
example is the Hopf manifold, diffeomorphic to $S^1\times S^{2n-1}$.

On a Vaisman manifold, the vector field
$\theta^\sharp-\sqrt{-1}J\theta^\sharp$ generates a one-dimensional
holomorphic, Riemannian, totally geodesic foliation. If this is regular
and if $M$ is compact, then the leaf space $B$ is a K\"ahler manifold. 

\hfill

\example
On a Hopf
manifold $\C^n\setminus\{0\}/\langle z_i\mapsto 2z_i\rangle$, the LCK
metric $\displaystyle\frac{\sum dz_i\otimes dz_i}{|\sum z_i\bar z_i|^2}$
is Vaisman and regular; the leaf space is $\C P^{n-1}$.

We refer to \cite{drag} or to the more recent \cite{ov11} for more details
about LCK geometry.

\hfill

It is known, \cite{_Tricerri_, vuli}, that the blow-up {\bf at points}
preserve the LCK class. The present paper is devoted to the blow-up of LCK
manifolds along subvarieties. In this case, the situation is a bit more
complicated and a discussion should be made according to the dimension of
the submanifold.

\hfill

\definition\label{defgck}
Let $Y\stackrel j\hookrightarrow M$ be a complex subvariety. We say
that $Y$ is {\bf of induced globally conformally K\"ahler type}
(IGCK) if the cohomology class $j^*[\theta]$ vanishes, where
$\theta$ denotes the cohomology class of the Lee form on $M$.

\hfill

\remark
Notice that a IGCK-submanifold of an LCK manifold is always K\"ahler.

\hfill

\remark\label{dim}
By a theorem of Vaisman (\cite{vaisman2}),
any LCK metric on a compact complex manifold $Y$
of K\"ahler type is globally conformally K\"ahler
if $\dim_\C Y>1$. Therefore, the IGCK condition
above for smooth $Y$ with $\dim_\C Y>1$ is
equivalent to $Y$ being K\"ahler.

\hfill

\remark\label{curve}
Notice that there may exist curves on LCK manifolds which are not IGCK, despite being obvioulsy of K\"ahler type. For instance, if $M$ is a regular Vaisman manifold, and if $Y$ is a fiber of its elliptic fibration, then $Y$ is not IGCK, as any compact complex subvariety of a compact Vaisman manifold has an induced Vaisman structure (see {\em e.g.} \cite[Proposition 6.5]{verbitsky_vanishing}).

\hfill

The main goal of the present paper is to prove the following
two theorems:

\hfill

\theorem\label{_main_blow-up_LCK_Theorem_}
Let $M$ be an LCK manifold, $Y\subset M$ be a
smooth complex {{IGCK subvariety}}, and let $\tilde M$
be the blow-up of $M$ centered in $Y$. Then
$\tilde M$ is LCK.

\hfill

{\bf Proof:} See the argument after \ref{chase}.
\endproof

\hfill

\theorem\label{_main_blow-down_Theorem_}
Let $M$ be a complex variety,
and $\tilde M \arrow M$ the blow-up of a compact
subvariety $Y\subset M$. Assume that
$\tilde M$ is smooth and admits an LCK metric. Then
the blow-up divisor $\tilde Y \subset \tilde M$
is a {IGCK subvariety}.

\hfill

{\bf Proof:} See \ref{_blow_up_Remark_}. \endproof

\hfill

\remark
In the situation described in \ref{_main_blow-down_Theorem_},
the variety $\tilde Y$ is of K\"ahler type, because it is IGCK.
When $Y$ is smooth, $Y$ is K\"ahler, as shown by Blanchard
(\cite[Th\'eor\`eme II.6]{blanchard}).
 {Together with \ref{dim},}
this implies the following corollary.

\hfill

\corollary\label{_blow-up_LCK_Corollary_}
Let $M$ be an LCK manifold, and $Y\subset M$ a smooth compact subvariety,
such that the blow-up of $M$ in $Y$ admits an LCK metric.
{If $\dim_\C(Y)>1$ } then $Y$ is a {IGCK subvariety}.
\endproof

\hfill


\remark\label{vaisman} Note that, from \cite[Proposition 6.5]{verbitsky_vanishing}, a compact
complex submanifold $Y$ of a compact Vaisman manifold is itself Vaisman,
and $\theta$ represents a non-trivial class in the cohomology of $Y$, 
{
so there are no IGCK submanifolds of proper dimension $\dim_\C(Y)>0.$}
This implies the following corollary.

\hfill

\corollary \label{_Vaisman_blow_Corollary_}
The blow-up of a compact Vaisman manifold along a compact
complex submanifold $Y$ of dimension at least {$1$} cannot have an LCK metric.

\hfill


The proofs of these two theorems and of the corollary will be
given in Section \ref{_Proofs_Section_}. As a
by-product of our proof, we obtain the following:

\hfill

\corollary\label{twistor}
If $M$ is a twistor space, and if  $M$ admits a LCK
metric, then this metric is actually GCK.

{\bf Proof:} See \ref{_twi_Corollary_}. \endproof

\hfill

Here, by a ``twistor space'' we understand any of 
the following constructions of a complex manifold: the twistor
spaces of half-conformally flat 4-dimensional Riemannian manifolds,
twistor spaces of quaternionic-K\"ahler manifolds, and Riemannian
twistor spaces of conformally flat manifolds.

\hfill

\remark $(i)$ A similar, weaker result is proven in
\cite{_Kokarev_Kotschick:Fibrations_}. Namely, the twistor
space of half-conformally flat 4-dimensional Riemannian manifolds with
large fundamental group cannot admit LCK metrics with automorphic
potential on the covering. The proof uses different techniques from ours,
and which cannot be generalized neither
to higher dimensions nor to quaternionic K\"ahler manifolds.

$(ii)$ It was known from \cite{paul, oleg} that the natural metrics (with
respect to the twistor submersion) cannot be LCK. Our result refers to any
metric on the twistor space, not necessarily related to the twistor
submersion.  On the other hand, as shown by Hitchin, the twistor space
of a compact 4-dimensional manifold is not of
K\"ahler type, unless it is
biholomorphic to $\C P^3$ or to the flag variety $F_2$ \cite{Hit}.

\hfill

\remark\label{rem1} 
So far, we were unable to deal with the reverse
statement of
\ref{_main_blow-up_LCK_Theorem_}, namely, to determine whether
a smooth bimeromorphic
contraction of an LCK manifold is always LCK.
In the particular case when 
an exceptional divisor is contracted 
to a point, this has been proven to be true by
Tricerri, \cite{_Tricerri_}; we
conjecture that in the general case this is 
false, but we are not able to find any
example.

For GCK (that is, K\"ahler) 
manifolds, the answer is well known: blow-downs of
K\"ahler manifolds can be non-K\"ahler, as one can see from any
example of a Moishezon manifold.

\hfill

\remark\label{rem2} 
{
We summarize the case of blow-up of curves on LCK manifolds. Since rational curves are simply-connected, they are IGCK submanifolds, so blowing-up a rational curve on a LCK manifold always yields a manifold of LCK type. The case of the elliptic curves was partially tackled in \ref{_Vaisman_blow_Corollary_}. If $Y$ is a curve of arbitrary genus contained in an exceptional divisor of a blow-up, then it is also automatically a IGCK subvariety  since the exceptional divisor  is so; hence again, blowing it up yields a manifold of LCK type.} 

To our present knowledge, the only examples of curves
$Y$ on LCK
manifolds $M$ with genus $g(Y)\geq 2$ are curves belonging to some
exceptional divisors.
It would be interesting to prove that this is the case in general, or to
build out a
counter-example.


\section{The proofs}
\label{_Proofs_Section_}


\lemma\label{fibr}
Let $M$ be an LCK manifold, $B$ a path connected topological space and let
$\pi:M\ra B$ be a continuous map. Assume that either
\begin{description}
\item[(i)] $B$ is an irreducible
complex variety, and $\pi$ is proper and holomorphic.
\item[(ii)] $\pi$ is a locally trivial fibration with fibers which
are complex subvarieties of $M$.
\end{description}
Suppose also that the map
$$\pi^*:H^1(B)\ra H^1(M)$$
 is an isomorphism, and the generic fibers
of $\pi$ are positive-dimensional.
Then the LCK structure on $M$ is actually GCK.

\hfill

{\bf Proof:}
Denote by $\theta$ the Lee form of $M$, and let
$\tilde M$ be the minimal GCK covering of $X$,
that is, the minimal covering $\tilde M \arrow M$ such that
the pullback of  $\theta$ is exact.
Since $H^1(B)\cong H^1(M)$, there exists a covering $\tilde B \arrow B$
such that the following diagram is commutative, {and the  fibers  of $\tilde \pi$ are compact}:
\[
\begin{CD}
\tilde M @>>> M\\
@V{\tilde \pi}VV @VV{\pi}V\\
\tilde B @>>> B
\end{CD}
\]
Let $\tilde B_0\subset \tilde B$
be the set of regular values of $\tilde \pi$,
and let $F_b:= \tilde \pi^{-1}(b)$ be the
regular fibers of $\tilde \pi$, $\dim_\C F_b =k$. Since
$B_0$ is connected, all $F_b$ represent
the same homology class in $H_{2k}(\tilde M)$.

Denote the K\"ahler form of $\tilde M$ by $\tilde \omega$,
conformally equivalent to the pullback of the Hermitian
form on $X$.

Since all $F_b$ represent the same homology class, the
Riemannian volume
\[
\Vol_{\tilde \omega}(F_b) := \int_{F_b}\tilde\omega^k
\]
is independent from $b\in B_0$. This gives (recall the definition of the
character $\chi$ in \eqref{chi})
\[
\Vol_{\tilde \omega}(F_b)= \int_{F_b}\tilde\omega^k
= \int_{F_{\gamma^{-1}(b)}}\gamma^*\tilde\omega^k =
\int_{F_{\gamma^{-1}(b)}}\chi(\gamma)^k \tilde\omega^k=
\chi(\gamma)^k\Vol_{\tilde \omega}(F_b),
\]
hence the constant $\chi_\gamma$ is equal to 1 for
all $\gamma \in \Gamma$. Therefore, $\tilde\omega$
is $\Gamma$-invariant, and $M$ is globally conformally
K\"ahler.
\endproof

\hfill

The above lemma immediately implies \ref{twistor}.

\hfill

\corollary\label{_twi_Corollary_}
Let $Z$ be the twistor space of $M$,
understood in the sense of \ref{twistor}.
Assume that $Z$ admits an LCK
metric. Then this metric  is globally conformally
K\"ahler.

\hfill

{\bf Proof:} There is a locally trivial fibration
$Z \arrow M$, with complex analytic fibers which are
compact symmetric K\"ahler spaces, hence
\ref{fibr} can be applied.
\endproof

\hfill

\remark\label{_blow_up_Remark_}
In the same way one deals with the blow-ups:
the generic fibers over an exceptional set 
of a blow-up map are positive-dimensional. Therefore,
\ref{fibr} implies \ref{_main_blow-down_Theorem_}.

\hfill

We can now give {\bf The proof of \ref{_Vaisman_blow_Corollary_}:} 

If $\dim_\C(Y)>1$ the result follows from  \ref{_blow-up_LCK_Corollary_} and \ref{vaisman}.
In the case $\dim_\C(Y)=1$ we cannot use this argument directly - see \ref{curve} - so in this case we  argue as follows.

 Assume $\tilde{M}$ has an LCK metric $\tilde \omega$ with Lee form $\tilde \eta$. By \ref{_main_blow-down_Theorem_}, the restriction $\tilde \eta_{|Z}$ to the exceptional divisor $Z$ is exact. Hence, after possibly making a conformal change of the LCK metric, we can assume $\tilde \eta_{|V}=0$ where $V$ is a neighbourhhood of $Z$. In particular, $\tilde\eta$ will be the pull-back of a one-form $\eta$ on $M$. On the other hand, $\tilde\omega$ gives rise to a current on $\tilde M$ (see also \S \ref{cur}) and its push-forward defines an LCK positive $(1,1)$ current $\Xi$ on $M$ with associate Lee form $\eta$. Clearly  $\eta_{|Y}=0$.
 
 Possibly conformally changing now $\Xi$, we can assume that $\eta$ is the unique harmonic form (with respect to the Vaisman metric of $M$) in its cohomology class.  Possibly $\eta_{|Y}$ is no longer zero, but remains {\em exact}. 
 
 We now show that $\eta$ is basic with respect to the canonical foliation $\mathcal{F}$ generated on $M$ by $\theta^\sharp-\sqrt{-1}J\theta^\sharp$. Indeed, from \cite{vaisman2}, we know that 
any harmonic form on a compact Vaisman manifold decomposes as a sum $\alpha+\theta\wedge\beta$ where $\alpha$ and $\beta$ are basic and transversally (with respect to $\mathcal{F}$) harmonic forms. In particular, as a transversally harmonic function is constant, we have
\begin{equation}\label{dec_har}
\eta=\alpha+ c\cdot\theta,
\end{equation}
where $c\in\RR$ and $\alpha$ is basic, transversally harmonic (see \cite{tondeur} for the theory of basic Laplacian and basic cohomology etc.). 

Let now $S^1$ denote the unique homology class in $H_1(M)$ (call it {\em the fundamental circle of $\theta$}) such that $\int_{S^1}\theta=1$ and $\int_{S^1}\alpha=0$ for every basic cohomology class $\alpha$. 

As any complex submanifold of a compact Vaisman manifold is tangent to the Lee field and hence Vaisman itself, $Y$ is Vaisman with Lee form $\theta_{|Y}$.  Hence  we deduce that the fundamental circle of $\theta$ is the image of the fundamental circle of $\theta_{\vert Y}$ under the natural map $H_1(Y)\rightarrow H_1(M)$. 

We now integrate \ref{dec_har} on any $\gamma\in H_1(Y)$ and take into account that $\eta_{|Y}$ is exact to get $c=0$. Hence, $\eta$ basic. It can then be treated as a harmonic one-form on a K\"ahler manifold (or use the existence of a transversal $dd^c$-lemma). This implies $d^c\eta=0$.
 
  But then one obtains a contradiction, as follows. Letting $J$ to be the almost complex structure of $M$, we see on one hand we have
$$\int_M d(\Xi^{n-1}) \wedge J(\theta)=\int_M (n-1)\Xi^{n-1} \wedge \theta\wedge J(\theta)>0$$
since $\Xi$ is positive. On the other hand, since $d(J(\theta))=0,$ it follows that $d(\Xi^{n-1}) \wedge J(\theta)$ is exact so $\int_M d(\Xi^{n-1}) \wedge J(\theta)=0,$ a contradiction.

\endproof

\hfill

The following result is certainly well-known, but since we were not able
to find out
an exact reference we include a proof here.

\hfill

\lemma\label{chase}
Assume $(U, g)$ is a K\"ahler complex manifold,
$Y\subset U$ a
compact submanifold and let $c:\tilde{U}\ra U$ be the blow-up of $U$ along
$Y$.
Then,  for any open neighbourhood $V\supset Y$, there is a K\"ahler metric
$\tilde{g}$ on
$\tilde{U}$ such that
$$\tilde{g}_{\vert \tilde{U}\setminus c^{-1}(V)}=c^*(g_{\vert U\setminus
V})$$

{\bf Proof.} (due to M. P\u aun; see also \cite{vuli}).

1. There is a (non-singular) metric on ${\mathcal O}_{\tilde{U}}(-D)$
(where $D$ is
the exceptional divisor of the blow-up) such that:

1.A. Its curvature is zero outside $c^{-1}(V)$, and

1.B. Its curvature is strictly positive at every point of $D$ and in any
direction tangent to $D$.

Indeed, if such a metric is found, everything follows, as the curvature of
this metric plus
a sufficiently large multiple of $c^*(g)$ will be positive definite on
$\tilde{U}.$

2. To finish the proof, we notice that the existence of a metric $h$  with
property 1.B
is clear, due to the restriction of ${\mathcal O}_{\tilde{U}}(-D)$ to $D$.

Now let
$\alpha$ be its curvature; then
$\alpha- i\ddbar \tau= -[D]$
for some function $\tau$,
with at most logarithmic poles along $D$, bounded from above, and
non-singular
on  $\tilde U\setminus D.$
Consider the function
$\tau_0:= \max(\tau, -C)$
where $C$
is some positive constant, big enough such that
on  $\tilde U\setminus c^{-1}(V)$ we have $ \tau> -C.$
Clearly, on a (possibly smaller) neighbourhood of $D$ we will have
$\tau_0= -C,$
such that the new metric
$e^{-\tau_0}h$
on ${\mathcal O}_{\tilde{U}}( -D)$
also satisfies 1.A.
\endproof

\hfill

{Now we can prove \ref{_main_blow-up_LCK_Theorem_}.
Let $c:\tilde M\ra M$ be the blow up of $M$ along the submanifold $Y$. Let
$g$ be a LCK metric on $M$ and let $\theta$ be its Lee form. 
Since $Y$ is IGCK we see
$\theta_{|Y}$ is exact.  Let
$U$ be a neighbourhood of $Y$ such that the inclusion $Y\hookrightarrow U$
induces an isomorphism of the first cohomology. 
Then $\theta_{\vert U}$ is also exact, 
so, after possibly conformally rescaling $g$, we may assume $\theta_{|U}=0$ and hence  $g_{\vert U}$
is K\"ahler. In particular,
$\mathrm{supp}(\theta)\cap U=\emptyset.$
Now choose a smaller neighbourhood $V$ of $Y$ and apply \ref{chase}. We
get a K\"ahler
metric $\tilde{g}$ on $\tilde{U}$
which equals $c^*(g)$ outside $c^{-1}(V)$, so it glues to $c^*(g)$ giving
a LCK
metric on $\tilde M$.
}
\endproof

\hfill

\noindent{\bf Acknowledgments.} We are indebted to M. Aprodu,
J.-P. Demailly, D. Popovici and M. Toma for useful discussions. We are
grateful to Mihai P\u aun for explaining us the proof of \ref{chase}.

{\small

\noindent {\sc Liviu Ornea\\
University of Bucharest, Faculty of Mathematics, \\14
Academiei str., 70109 Bucharest, Romania. \emph{and}\\
Institute of Mathematics ``Simion Stoilow" of the Romanian Academy,\\
21, Calea Grivitei Street
010702-Bucharest, Romania }\\
\tt Liviu.Ornea@imar.ro, \ \ lornea@gta.math.unibuc.ro

\hfill

\noindent {\sc Misha Verbitsky\\
{\sc  Laboratory of Algebraic Geometry,
Faculty of Mathematics, NRU HSE,
7 Vavilova Str. Moscow, Russia }\\
\tt verbit@maths.gla.ac.uk, \ \  verbit@mccme.ru

\hfill

\noindent {\sc Victor Vuletescu\\ University of Bucharest, Faculty of
Mathematics,
\\14
Academiei str., 70109 Bucharest, Romania.}\\
\tt vuli@gta.math.unibuc.ro
}
}

\end{document}